\documentstyle[12pt]{article}

\input amssym.def
\input amssym

\newtheorem{prop}{Proposition}
\newtheorem{cor}{Corollary}
\newtheorem{trm}{Theorem}

\newcommand\OSup{\bigoplus\limits}

\newcommand\<{\preccurlyeq }

\def\>{\succcurlyeq }

\def\0{{\bf 0}}
\def\1{{\bf 1}}

\def\a-{$a\mbox -$\nolinebreak }
\def\aa-{$a\mbox -$\nolinebreak }
\def\b-{$b\mbox -$\nolinebreak }
\def\c-{$\vee\mbox -$\nolinebreak }

\def\d-{$\delta\mbox -$\nolinebreak }

\def\oo-{$o\mbox -$\nolinebreak }
\def\wo-{$wo\mbox -$\nolinebreak }
\def\G-{$G\mbox -$\nolinebreak }
\def\g-{$g\mbox -$\nolinebreak }

\def\R-{{${\bf R_{\oplus }}\mbox -$\nolinebreak }}
\def\plus-{{$\oplus\mbox -$\nolinebreak}}

\def\Rmax{\Bbb R_{\max}}

\def\Rmin{\Bbb R_{\min}}
\def\Zmax{\Bbb Z_{\max}}

\def\bC{\Bbb C}
\def\Conc{\mbox{Conc}}
\def\Mat{\mbox{\rm Mat}}
\newcommand\Lip{{\rm Lip}}
\newcommand\lip{{\rm lip}}

\begin{document}

\begin{center}
{\large{\bf KERNEL THEOREMS AND NUCLEARITY IN
IDEMPOTENT MATHEMATICS.\\[2mm] 
AN ALGEBRAIC APPROACH}}\\[3mm]

G.L.~Litvinov and G.B.~Shpiz
{\footnote{This work has been supported by the RFBR-CNRS grant 05-01-02807.\\
The paper will be published in Journal of Mathematical Sciences (New York).}}
\\[2mm]

\end{center}

{\qquad\qquad\qquad\qquad\qquad\qquad Dedicated
 in dear memory of F.A.~Berezin}

\begin{abstract} In the framework of idempotent mathematics,
analogs of the classical kernel theorems of L. Schwartz and
A.~Grothendieck are studied. Idempotent versions of nuclear
spaces (in the sense of A.~Grothendieck) are discussed. The
so-called algebraic approach is used. This means that the
basic concepts and results (including those of ``topological'' 
nature) are simulated in purely algebraic terms.
\end{abstract}

2000{\it Mathematics Subject Classification.} Primary 46S10,
46T99, 47H99, 06A99, 06F99; Secondary 06B35, 06F05.


\section*{Introduction}

Idempotent mathematics is based on replacing the usual 
numerical fields with idempotent semifields and semirings.
In other words, the usual arithmetic operations are replaced by
a new set of basic associative operations (a new addition $\oplus$ 
and a new multiplication $\odot$) and all the semifield axioms
or semiring axioms hold; moreover the new addition is
idempotent, i.e. $x\oplus x=x$ for every element $x$ of the
corresponding semiring, see, e.g., \cite{1}--\cite{10}. A typical
example is the semifield $\Rmax=\Bbb R\bigcup\{-\infty\}$
known as the Max-Plus algebra. This semifield consists of all real
numbers and an additional element $\0=-\infty$.
This element $\0$ is the zero element in $\Rmax$ and the basic
operations are defined by the formulas $x\oplus y=\max\{x,y\}$ 
and $x\odot y=x+y$; the identity (or unit)
element~$\1$ coincides with the
usual zero $0$. Similarly the semifield $\Rmin=\Bbb R\bigcup\{+\infty\}$
is defined. It consists of all real numbers and an additional element
$\0=+\infty$; the basic operations are $\oplus=\min$ and $\odot=+$.
Of course, the semifields $\Rmax$ and $\Rmin$ are isomorphic. Many
nontrivial examples of idempotent semirings and semifields can 
be found, e.g., in \cite{3}--\cite{10}.

Linear algebra over idempotent semirings was constructed
by many authors starting from S.C.~Klenee, N.N.~Vorobjev,
B.A.~Carre, R.A.~Cuningham-Green and others. Basic concepts
and results of idempotent analysis and functional analysis are
established 
by V.P.~Maslov and his collaborators, see, e.g., 
\cite{1}--\cite{14}; in particular, important results and applications of
idempotent analysis are due to V.N.~Kolokoltsov (see, e.g.,
\cite{3}--\cite{6}, \cite{9}, \cite{13}--\cite{15}). 
An important development of these
ideas is due to the 
French mathematicians M.~Akian, G.~Cohen,
S.~Gaubert, J.-P.~Quadrat and others, see, e.g., 
\cite{9, 15, 16} and the survey \cite{10}. 
A remarkable paper of O.~Viro \cite{17} was a starting point
for idempotent 
geometry; now this subject is usually called
tropical geometry. The subject is very popular and substantial
contributions to its development are due to many authors
(M.~Kontsevich, G.~Mikhalkin, M.~Kapranov, I.~Itenberg,
V.~Kharlamov, E.~Shustin, B.~Sturmfels and others).
There are other areas of idempotent mathematics, e.g.,
idempotent interval analysis \cite{18, 19}. Concerning
the further development of the subject and its history,
see, e.g., the survey~\cite{10}.

Mathematics over idempotent semirings can 
be treated as
a result of a dequantization of the traditional mathematics
(over fields) and as 
its peculiar ``shadow.'' This ``shadow''
stands in the same relation to the traditional mathematics
as does classical physics to quantum theory (see details in
\cite{7, 8, 10}). There is a heuristic correspondence between
important, interesting, and useful constructions and results of the
traditional mathematics over fields and similar constructions
and results over idempotent semifields and semirings \cite{7}. 
In many respects idempotent mathematics is simpler than the
traditional one. However, the transition from traditional
concepts and results to their idempotent versions is often
nontrivial.

The aim of the paper is to describe idempotent versions of the 
classical kernel theorems of L.~Schwarts and A.~Grothendieck 
(see, e.g., \cite{21, 22}). In other words, we describe conditions
under which, in idempotent functional analysis, linear operators have
integral representations in terms of the 
idempotent integral of
V.P.~Maslov. We define the notion of nuclear idempotent semimodule
similar to the notion of nuclear space (in the sense of
A.~Grothendieck) in the traditional analysis (see, e.g., 
\cite{21, 22}). Moreover, we give a 
rather explicit description
of semimodules of functions for which our kernel theorem
is true. We use the so-called algebraic approach; this means that
the basic concepts and results (including those of topological
nature) are simulated in purely algebraic terms; in our 
subsequent papers a topological approach will be also used.
This paper continues the series of publications on idempotent
functional analysis \cite{8, 11, 12} and we use the notation and terminology
defined in those articles. Some of our results on idempotent
semimodules of functions depend on 
semimodules' imbeddings
into the corresponding function spaces. In this case sometimes
we discuss nonlinear (in the sense of idempotent mathematics)
mappings and functionals. There are other results which are
invariant with respect to such embeddings and can be rewritten
in an invariant form.

For some concrete idempotent semimodules consisting of
continuous or bounded functions, concrete kernel
theorems are presented in \cite{1}--\cite{5}, \cite{13}, 
\cite{23}--\cite{25}. 
A limiting case for kernel theorems is a result on integral
representations of linear functionals; results of this
kind are presented in \cite{1}--\cite{6}, \cite{8, 12, 13}, 
\cite{23}--\cite{30}. 
A generalization of these results and their unification will
be examined in a separate paper. In \cite{31}, for the idempotent
semimodule of all bounded functions with values in $\Rmin$,
there was posed a problem of describing the class of
subsemimodules where the corresponding kernel theorem holds as well
as of 
describing the corresponding integral representations of linear
operators. In the present paper a very general case
of semimodules over boundedly complete idempotent semirings is
examined. Some of the results obtained here can be regarded
as possible versions of an answer to the problem posed in
\cite{31}. Basic results of this paper were announced in~\cite{32}.

\section{Idempotent semimodules}

\subsection{Notation and basic terminology}
In the present paper we use terminology and notation from
\cite{8, 11, 12}. Recall that an {\it idempotent semigroup}
is an additive semigroup with commutative addition $\oplus $ 
such that for every element $x$ the equality $x\oplus x=x$ holds.
Every idempotent semigroup will be treated as an ordered  set
with respect to the following (partial) {\it standard order}:
$x\< y$ if and only if $x\oplus y=y$. It is easy to see that this
order is well-defined and $x\oplus y=\sup \{x,y\}$,
that is every idempotent semigroup is an upper semilattice
{\footnote{Let us remind that a partially ordered set $V$
is called an {\it upper semilattice} if for arbitrary elements $x$, 
$y$ of $V$ there exists their least upper bound $\sup\{x,y\}$
in $V$.}} with
respect to the standard order (see, e.g., \cite{33}).
For an arbitrary subset $X$ of an idempotent semigroup we set
$\oplus X=\sup X$ and $\wedge X=\inf X$ if $\sup X$ and $\inf X$
exist; in particular, we suppose that $X$ is bounded from
above or below respectively. An idempotent semigroup is called
{\it \b-complete} (or {\it boundedly complete}) if every its subset
that is 
bounded from above (including the empty set) has the least upper
bound. In particular, every \b-complete idempotent semigroup
contains a zero element (denoted $\0$) which coincides with
$\oplus \varnothing$, where $\varnothing$ is the empty set.
Thus every nonempty subset of this semigroup is bounded
from below by the zero $\0$. So every nonempty subset of the \b-complete
semigroup is bounded if and only if this subset is
bounded from above. A homomorphism $g$ for \b-complete
idempotent semigroups is called a {\it\b-homomorphism}, if
 $g(\oplus X)=\oplus g(X)$ for every subset $X$ bounded 
from above.
{\it An idempotent semiring} is an idempotent semigroup 
endowed with an associative multiplication $\odot $ with
an identity element (denoted $\1 $) such that the
corresponding distribution laws (left and right) hold.
An idempotent semiring $K$ is {\it \b-complete} if it is a
\b-complete idempotent semigroup and the following ``infinite''
distribution laws hold: $k\odot \oplus X=\oplus (k\odot X)$ 
and $(\oplus X)\odot k=\oplus (X\odot k)$ for every
$k\in K$ and every bounded subset $X$ of $K$.

A commutative idempotent semiring is called an {\it idempotent
semifield} if every its nonzero element is invertible. A
semifield is a \b-complete idempotent semiring if and only if it
is a \b-complete idempotent semigroup
with respect to the operation~$\oplus$~\cite{8}.
The algebra
 $\Rmax$ (or $\Rmin$) described in the Introduction is an
example of a \b-complete semifield. Another example is the
semifield $\Zmax$ consisting of integer elements of $\Rmax$ 
and the element $\0=-\infty$ with operations induced from
 $\Rmax$.

If we add a maximal 
element $+\infty$ to $\Rmax$ and extend the operations 
in an obvious way, 
we shall get a new idempotent semiring (not a semifield!)
 $\widehat{\Bbb R}_{\max} = \Rmax\cup \{+\infty\}$;
the semiring $\widehat{\Bbb Z}_{\max} = \Zmax\cup \{+\infty\}$ is
defined similarly. These semirings  
are not only \b-complete but complete
as partially ordered sets (with respect to the standard order
$\< $) and the passage from $\Rmax$ and $\Zmax$ to $\widehat{\Bbb R}_{\max}$
and $\widehat{\Bbb Z}_{\max}$ is the so-called normal completion in
the sense of \cite{33}.  Square matrices of the size $n\times n$
with entries from a \b-complete idempotent semiring $K$ form a
\b-complete noncommutative semiring $\Mat_n (K)$ with respect
to the operations of addition and multiplication of matrices
generated by the operations in $K$.

An {\it idempotent semimodule} over an idempotent semiring $K$
is an additive commutative idempotent semigroup $V$, with the
addition operation denoted by $\oplus$ and the zero element
denoted by $\0$, such that a (left) multiplication $k\odot x$
is defined for all $k\in K$ and $x\in V$ 
in such a way that
the usual rules are satisfied:
$a\odot (b\odot x)=(a\odot b)\odot x$,
$(a\oplus b)\odot x=a\odot x\oplus b\odot x$,
$a\odot (x\oplus y)=a\odot x\oplus a\odot y$,
$\0 \odot x=\0$, $\1\odot x = x$ for all $a,b \in K$ and $x,y \in V$. 

In what follows we shall suppose that all the semigroups, semirings,
semifields, and semimodules are idempotent unless
otherwise specified.

A semimodule $V$ over a \b-complete semiring $K$ is {\it \b-complete}
if it is a \b-complete semigroup and the following infinite
distribution laws hold:
$k\odot \oplus X=\oplus (k\odot X)$ for every $k\in K$ and every
bounded subset $X\subset V$ as well as $(\oplus Q) \odot v=\oplus
(Q\odot v)$ for every bounded subset $Q$ of $K$ and every element
$v\in V$ (see \cite{8}, definition 4.3).

A homomorphism $g: V\to W$ for \b-complete semimodules $V$ and $W$
is called a {\it \b-homomorphism} or a {\it \b-linear mapping}
({\it operator}) if $g(\oplus X)=\oplus g(X)$ for every bounded
subset $X\subset V$ (a more general definition for the case of
incomplete semimodules see in \cite{8}); of course, here
$g(a\odot x)= a\odot g(x)$ for all $a\in K$, $x\in V$. 
Homomorphisms taking their values in the basic semiring $K$
(treated as a semimodule over itself) are called 
{\it linear functionals}. Of course, \b-{\it linear functionals}
are linear functionals such that the corresponding homomorphisms are
\b-homomorphisms.

In idempotent mathematics,  \b-linear operators and functionals can
be regarded as analogs of traditional (semi)continuous linear
operators and functionals in classical analysis, see~\cite{8}.

Let $V$ and $W$ be \b-complete semimodules over a \b-complete
semiring $K$. Denote by $L_b(V,W)$ the set of all \b-linear
mappings from $V$ to $W$. It is easy to check that $L_b(V,W)$ 
is an idempotent semigroup with respect to the pointwise
addition of operators; the composition (product) of \b-linear
operators is also a \b-linear operator and therefore 
the set $L_b(V,V)$ 
is an idempotent semiring with respect to these operations,
see, e.g., \cite{8}. The following proposition can be treated
as a version of the Banach-Steinhaus theorem in idempotent
analysis.

\begin{prop}Suppose that $S$ is a subset in $L_b(V,W)$ and the set
$\{g(v)\mid g\in S\}$ is bounded in $W$ for every element
$v\in V$; thus the element $f(v) = \OSup_{g\in S} g(v)$ 
exists because the semimodule $W$ is \b-complete. Then the
mapping $v\mapsto f(v)$ is a \b-linear operator, i.e. an
element of $L_b(V,W)$. The subset $S$ is bounded; moreover, $\sup S = f$.
\end{prop}

{\bf Proof.} It is easy to check that $f$ is a homomorphism $V\to W$.
For every bounded subset $V$ in $W$ the following equations hold:
$f(\oplus X) = \OSup_{g\in S} g(\oplus X) = \OSup_{g\in S} \OSup_{v\in X}
g(v) = \OSup_{v\in X} \OSup_{g\in S} g(v) = \OSup_{v\in X} f(v) =
\oplus f(X)$. This means that $f$ is a \b-linear operator. 
By our construction, $f = \sup S$ in $L_b(V,W)$.$\;\square$

\begin{cor}The set $L_b(V,W)$ is a \b-complete idempotent semigroup
with respect to the (idempotent) pointwise addition of
operators. If $V = W$, then $L_b(V,V)$ is a \b-complete idempotent
semiring with respect to the operations of pointwise addition and
composition of operators.
\end{cor}

\begin{cor} A subset $S$ is bounded in $L_b(V,W)$  if and only if
the set $\{g(v)\mid g\in S\}$ is bounded in the semimodule $W$
for every element $v\in V$.
\end{cor}

These corollaries can be easily deduced from Proposition 1 and
the basic definitions.

A subset of an idempotent semimodule is called a {\it subsemimodule}
if it is closed under addition and multiplication by scalar coefficients.
A subsemimodule $V$ of a \b-complete semimodule $W$ is {\it \b-closed}
if $V$ is closed with respect to summing of all subsets of $V$ 
that are bounded in $W$. 
A subsemimodule of a \b-complete semimodule is called a
{\it \b-subsemimodule} if the corresponding embedding is a \b-homomorphism.
It is easy to see that each \b-closed subsemimodule is a \b-subsemimodule
but the converse is not true (see Section~1.2). 
The main feature of \b-subsemimodules is that restrictions of \b-linear
operators and functionals to these semimodules are \b-linear. 

The following definitions are very important for our aims.
Suppose that $W$ is an idempotent \b-complete semimodule over a
\b-complete idempotent semiring $K$ and $V$ is a subset of $W$ such
that $V$ is closed under multiplication by scalar coefficients
and is an upper semilattice with respect to the order induced
from $W$. Let us define an addition operation in $V$ by the
formula $x\oplus y = \sup\{ x, y \}$, where $\sup$ means the least 
upper bound in $V$. If $K$ is a semifield, then $V$ is a
semimodule over $K$ with respect to this addition.

For an arbitrary \b-complete semiring $K$ we shall say that
$V$ is a {\it quasisubsemimodule} of $W$ if $V$ is a
semimodule with respect to this addition (this means that the
corresponding distribution laws hold).

A quasisubsemimodule $V$ of an idempotent \b-complete semimodule
$W$ is called a $\wedge$-{\it subsemimodule} if it contains 
$\0$ and is closed under the operations of taking infima (greatest lower
bounds) in $W$. It is easy to check (e.g., using lemma 2.1 in
\cite{8}), that {\it each $\wedge$-subsemimodule is a \b-complete 
semimodule}.

Note that quasisubsemimodules and $\wedge$-subsemimodules 
may fail to be subsemimodules because only the order is induced
and not the
corresponding addition (see Example 6 below).

Following \cite{8}, we say that idempotent semimodules over semirings
are {\it idempotent spaces}. In idempotent mathematics, such spaces
are analogs of traditional linear (vector) spaces over fields.
In a similar way we use the corresponding terms like
{\it \b-spaces, \b-subspaces, \b-closed subspaces}, $\wedge$-{\it subspaces}
etc.

Numerous examples of idempotent semimodules and spaces can be
found in \cite{8}; see also \cite{3}--\cite{9}, \cite{11}, 
\cite{14}--\cite{16}, \cite{25}. Some examples
are presented below, see, e.g., Section 3. 

\subsection{Functional semimodules}
Let $X$ be an arbitrary nonempty set and $K$ be an idempotent semiring.
By $K(X)$ denote the semimodule of all mappings (functions) $X \to K$
endowed with the pointwise operations. By $K_b(X)$ denote the subsemimodule
of $K(X)$ consisting of  all bounded
mappings. If $K$ is a \b-complete semiring, then $K(X)$ and $K_b(X)$ are
\b-complete semimodules. Note that $K_b(X)$ is a \b-subsemimodule but not
a \b-closed subsemimodule of $K(X)$. For any point $x\in X$, by
$\delta _x$ denote the functional on $K(X)$ that maps 
$f$ to $f(x)$. It
can be easily checked that the functional $\delta _x$ is \b-linear on $K(X)$.

We shall say that any quasisubsemimodule of $K(X)$ is an (idempotent)
{\it functional semimodule} on the set $X$. An idempotent functional
semimodule in $K(X)$ is called {\it \b-complete} if it is a \b-complete
semimodule.

A functional semimodule $V\subset K(X)$ is called a {\it functional
\b-semimodule} if it is a \b-subsemimodule of $K(X)$; a functional
semimodule $V\subset K(X)$ is called a {\it functional $\wedge$-semimodule}
if it is a $\wedge$-subsemimodule of $K(X)$.

In general, a functional of the form $\delta _x$ on a functional semimodule 
is not even linear, much less \b-linear 
(see Example~6 below). However, the following
proposition holds, which is a direct consequence of our definitions.

\begin{prop} An arbitrary \b-complete functional semimodule 
$W$ on a set $X$ is a \b-subsemimodule of $K(X)$ if and only if
each functional of the form $\delta _x$ (where $x\in X$) is
\b-linear on $W$.\end{prop}

{\bf Proof.}  Both statements mean that for each set of
functions $S\subset W$ bounded in $W$  the formula
$(\oplus S)(x)= (\oplus S(x))$ holds for all $x\in X$. $\;\square$\\

{\bf Example 1.} The semimodule $K_b(X)$ (consisting of all bounded
mappings from an arbitrary set $X$ to a 
\b-complete idempotent semiring $K$)
is a functional $\wedge$-semimodule. Hence it is a \b-complete
semimodule over $K$. Moreover, $K_b(X)$ is a \b-subsemimodule of the
semimodule $K(X)$ consisting of all mappings $X\to K$.\\

{\bf Example 2.} If $X$ is a finite set consisting of $n$ elements
($n>0$), then $K_b(X) = K(X)$ is an ``$n$-dimensional'' semimodule
over $K$ and it is denoted by $K^n$. In particular, $\Rmax^n$ is an
idempotent space over the semifield $\Rmax$ and $\widehat{\Bbb R}_{\max}^n$ 
is a semimodule over the semiring $\widehat{\Bbb R}_{\max}$. Note
that $\widehat{\Bbb R}_{\max}^n$ can be treated as a space
over the semifield $\Rmax$. For example, the semiring
 $\widehat{\Bbb R}_{\max}$ can be treated as a space (semimodule) over
$\Rmax$.\\

{\bf Example 3.} Let $X$ be a topological space. Denote by $USC(X)$ 
the set of all upper semicontinuous functions taking their
values in $\Rmax$. By definition, a function $f(x)$ is upper
semicontinuous if the set $X_s =\{ x\in X\mid f(x)\geq s\}$
is closed in $X$ for every element $s\in \Rmax$ (see, e.g.,  
\cite{8}, section 2.8). If a collection $\{f_{\alpha}\}$ consists of
upper semicontinuous (e.g., continuous) functions and  
$f(x) = \inf_{\alpha} f_{\alpha} (x)$, then $f(x) \in USC(X)$.  
It is easy to check that $USC(X)$  has a natural structure 
of an idempotent space over $\Rmax$. Moreover, $USC(X)$ is a functional
$\wedge$-space on $X$ and \b-space. The subspace
$USC(X)\cap K_b(X)$ of $USC(X)$ consisting of bounded (above) 
functions has the same properties.\\ 

{\bf Example 4.} Suppose that  $X$ is a partially ordered set and 
$K$ is the semifield $\Rmax$ (or $\Zmax$). Denote by $N(X)$ the set
of all nonincreasing functions defined on $X$ and taking their values
in $K$. It is easy to check that $N(X)$ has a natural structure
of a space over $K$; moreover, $N(X)$ is a functional \b-space
and a functional $\wedge$-space on $X$. The space $N(X)\cap K_b(X)$
has the same properties.\\

{\bf Example 5.} Let $\bC$ be the set of all complex numbers with
a singled out 
axis $\Bbb R$ of all real numbers. Denote by 
$N_{\Bbb R}(\bC)$ the set of all functions $\bC\to K$ nonincreasing
along the real axis $\Bbb R$, where $K = \Rmax$ or $K = \Zmax$.
It is easy to check that $N_{\Bbb R}(\bC)$ (as well as 
$N_{\Bbb R}(\bC)\cap K_b(\bC)$) has a natural structure of a space over 
 $K$ and this space is a functional \b-space and a functional
$\wedge$-space on $X$. This is an example of a natural ``intersection''
of space structures.\\

{\bf Example 6.} Note that an idempotent functional semimodule
(and even a functional $\wedge$-semimodule) on a set $X$ is not
necessarily a subsemimodule of $K(X)$. The simplest example is the
functional space (over $K=\Rmax$) $\Conc(\Bbb R)$ consisting of all 
concave functions on $\Bbb R$ with values in $\Rmax$.
Recall that a function $f$ belongs to $\Conc(\Bbb R)$ 
if and only if the subgraph
of this function is convex, i.e., the formula
$f(ax+(1-a)y)\ge af(x)+ (1-a)f(y)$ is valid for $0\le a\le 1$.
The basic operations with $\0\in \Rmax$ can be defined in an obvious way.
If $f,g \in \Conc(\Bbb R)$, then denote by $f\oplus g$ the sum of 
these functions in $\Conc(\Bbb R)$. The subgraph of $f\oplus g$ is the
convex hull of the subgraphs of $f$ and $g$. Thus $f\oplus g$
does not coincide with the pointwise sum (i.e., $\max\{f(x), g(x)\}$).

\section{Integral representations of linear operators}

\subsection{Integral representations of linear operators in\\
functional semimodules} 
Suppose $W$ is an idempotent \b-complete semimodule over a \b-complete
semiring $K$ and $V\subset K(X)$ is a \b-complete functional
semimodule on $X$. A mapping $A:V\to W$ is called an {\it integral operator}
or an operator with an {\it integral representation} if there exists
a mapping $k:X\to W$, called the {\it integral kernel} (or {\it kernel})
{\it of the operator} $A$, such that 
\begin{equation}
Af=\OSup_{x\in X} f(x)\odot k(x). 
\end{equation}
In idempotent analysis, the right hand side of formula (1) is often
written as $\int_X^{\oplus}f(x)\odot k(x) dx$ (see \cite{1}--\cite{7}).
Regarding the kernel $k$ it is supposed that the set 
$\{f(x)\odot k(x)|x\in X\}$ is bounded in $W$ for all $f\in V$ and $x\in X$.
We denote the set of all functions with this property by 
$kern_{V,W}(X)$. In particular, if $W=K$ and $A$ is a functional,
then this functional is called {\it integral}. Thus each integral
functional can be presented in the form of a ``scalar product''
$f \mapsto \int_X^{\oplus} f(x) \odot k(x)\; dx$, where 
$k(x)\in K(X)$; in idempotent analysis this situation is
standard, see, e.g.,~\cite{1}--\cite{10}.

Note that the functional of the form $\delta_y$ (where $y\in X$)
is a typical integral functional; in this case
$k(x) = \1$ if $x = y$ and $k(x) = \0$ otherwise.

We call a functional semimodule $V\subset K(X)$ {\it nondegenerate} 
if for each point $x\in X$ there exists a function $g\in V$
such that $g(x)=\1$, and {\it admissible} if for each function $f\in V$ 
and each point $x\in X$ such that $f(x)\neq \0$ there exists a 
function $g\in V$ such that $g(x)=\1$ and $f(x)\odot g\<f$.

Note that all idempotent functional semimodules over semifields
are admissible (it is sufficient to set $g = f(x)^{-1}\odot f$).

\begin{prop} Denote by $X_V$  the subset of $X$ defined by the
formula $X_V=\{ x\in X\mid\; \exists f\in V: f(x)=\1 \}$. If the semimodule 
$V$ is admissible, then the restriction to $X_V$ defines an
imbedding $i:V\to K(X_V)$ and its image $i(V)$ is admissible and
nondegenerate.

If a mapping $k:X\to W$ is a kernel of a mapping $A:V\to W$, then the
mapping $k_V:X\to W$ which is equal to $k$ on $X_V$ and equal to $\0$ 
on $X-X_V$ is also a kernel of the mapping $A$. 

A mapping $A:V\to W$ is integral if and only if the mapping
 $i_{-1}A:i(A)\to W$ is integral.
\end{prop}

{\bf Proof.} If the semimodule $V$ is admissible, then it is easy
to see that  $f(x) =\0$ for all $f\in V$ and $x\in X-X_V$. All the
statements of Proposition 3 can be easily deduced from this assertion.
$\;\square$\\

In what follows, $K$ always denotes a fixed \b-complete idempotent
semiring. We shall discuss semimodules over 
this semiring. If an operator has an integral representation,
this representation may not be unique. However, if the semimodule
$V$ is nondegenerate, then the set of all kernels of a fixed
integral operator is bounded with respect to the natural order
in the set of all kernels and is closed under the supremum
operation applied to its arbitrary subsets. In particular, 
{\it any integral
operator defined on a nondegenerate functional semimodule has
a unique maximal kernel}.

An important point is that an integral operator is not necessarily \b-linear
and even linear except when $V$ is a \b-subsemimodule of $K(X)$
(see Proposition 4 below).

If $W$ is a functional semimodule on a nonempty set $Y$, 
then the integral kernel $k$ of an operator $A$ can be naturally identified
with a function on $X\times Y$ defined by the formula 
$k(x,y)=(k(x))(y)$. This function will also be called the
{\it integral kernel} (or {\it kernel}) of the operator $A$.
As a result, the set  $kern_{V,W}(X)$  is identified with the set
$kern_{V,W}(X,Y)$ of all mappings $k: X\times Y\to K$ such that for every
point $x\in X$ the mapping $k_x: y\mapsto k(x,y)$ lies in $W$ and for
every $v\in V$ the set $\{v(x)\odot k_x|x\in X\}$ is bounded in $W$.
Accordingly, the set of all integral kernels of \b-linear
operators can be embedded to $kern_{V,W}(X,Y)$.

If $V$ and $W$ are functional \b-semimodules on $X$ and $Y$,
respectively, then the set of all kernels of \b-linear operators can
be identified with $kern_{V,W}(X,Y)$ (see Proposition 4 below)
and the following formula holds:
\begin{equation}
Af(y)=\OSup_{x\in X} f(x)\odot k(x,y)=\int_X^{\oplus} f(x)
\odot k(x,y) dx.
\end{equation}
This formula coincides with the usual definition of an
operator's integral representation. Note that formula (1) can
be rewritten in the form
\begin{equation}
Af=\OSup_{x\in X} \delta _x(f)\odot k(x). 
\end{equation}

\begin{prop} An arbitrary \b-complete functional semimodule
$V$ on a non\-empty set $X$ is a functional \b-semimodule on
$X$ (i.e., a \b-subsemimodule of $K(X)$) if and only if
all integral operators defined on $V$ are \b-linear.
\end{prop}

{\bf Proof.} Let $W$ be an arbitrary \b-complete semimodule. If $V$ 
is a functional \b-semimodule on $X$, then every operator 
$\Delta_{x,w}:V\to W$ of the form $f\mapsto \delta _x(f)\odot w$, 
where $x\in X$, $w\in W$ is \b-linear by virtue of Proposition 2. 
By definition each integral operator is a sum of operators of
this type; so every integral operator is \b-linear by virtue
of Proposition 1 and its corollaries. On the other hand, if each integral
operator defined on $V$ is \b-linear, then all the functionals
of the form $\delta _x$ are \b-linear because these functionals are
integral operators from $V$ into $K$. From this it follows that
$V$ is a functional \b-semimodule on $X$ because of Proposition 2.
So the proposition is proved.$\;\square $\\

The following concept (definition) is especially important
for our aims. Let $V\subset K(X)$ be a \b-complete functional
semimodule over a \b-complete idempotent semiring $K$. We
shall say that the {\it kernel theorem} holds for the
semimodule $V$ if each \b-linear mapping from $V$ into an 
arbitrary \b-complete semimodule over $K$ has an integral
representation. 

\begin{trm} Suppose that a \b-complete semimodule $W$ over a 
\b-complete semiring $K$ and an admissible functional $\wedge$-semimodule 
$V\subset K(X)$ are given.  Then each \b-linear operator
$A:V\to W$ has an integral representation of the form $(1)$. In particular,
if $W$ is a functional \b-semimodule on a set $Y$, then the
operator $A$ has an integral representation of the form $(2)$.
So for the semimodule $V$ the kernel theorem holds.
\end{trm}

{\bf Proof.} Denote by $X_V$ the subset of $X$ defined by
the formula $X_V = \{x\in X|\; \exists f\in V : f(x) =\1\}$
(see Proposition 3). If $x\in X_V$, then we
set $d_x=\wedge\{f\in V|f(x)=\1\}$. 
By our construction we have $f(x)\odot d_x\< f$. The semimodule $V$ is
admissible and it is a $\wedge$-semimodule, so $d_x(x)=\1$. 
From this we can easily deduce that $f=\OSup_{x\in X_V}f(x)\odot d_x$. 
Then $A(f)=\OSup_{x\in X_V}f(x)\odot A(d_x)$, that is the
mapping $x\mapsto k(x)$ 
defined by the formulas $k(x)=A(d_x)$ for all $x\in X_V$
and $k(x) = 0$ for each $x\notin X_V$
is a kernel of the operator $A$. The theorem is proved.$\quad\square$\\

{\bf Remark 1.} If in the framework of Theorem 1 the semimodule $V$ 
is nondegenerate, then the function $x\mapsto d_x$ is the maximal
integral kernel of the identity operator $id:V\to V$.

Indeed, if under the conditions of Theorem 1 the semimodule $V$ is
nondegenerate, then all the functions of the form $x\mapsto d_x$ 
belong to $V$. So if $k: X\to V$ is an integral kernel of the identity
operator $id:V\to V$, then $d_y=\OSup_{x\in X} d_y(x)\odot k(x)
\> d_y(y)\odot k(y)$ for each $y\in X$.  So we have $d_y\> k(y)$ for
each $y\in X$ because $d_y(y)=\1$, as was to be proved.\\

{\bf Remark 2.} Examples of admissible functional
$\wedge$-semimodules (and $\wedge$-spaces) appearing in Theorem 1
are presented above in the end of Section 1. Thus for these functional
semimodules and spaces $V$ over $K$ the kernel theorem holds
and each \b-linear operator mapping $V$ into an arbitrary \b-complete
semimodule $W$ over $K$ has an integral representation (2). 
Recall that each functional space over a \b-complete semifield is
admissible, see above.

\subsection{Integral representations of b-nuclear operators}
Let us introduce some important definitions. Suppose that $V$ and $W$
are \b-complete semimodules. A mapping $g: V\to W$ is called
{\it one-dimensional} (or a {\it mapping of the rank} 1) if it is of
the form $v\mapsto \phi (v)\odot w$, where $\phi$ is a \b-linear 
functional on $V$ and $w\in W$. A mapping $g$ is called {\it \b-nuclear}
if it is a sum of a bounded set of one-dimensional mappings. Each
one-dimensional mapping is \b-linear because the functional 
$\phi$ is \b-linear, so {\it every \b-nuclear operator is \b-linear}
(see Corollary 1 above). Of course, \b-nuclear mappings are closely
related with tensor products of idempotent semimodules, see~\cite{11}.

By $\phi\odot w$ we shall denote the one-dimensional operator
$v\mapsto\phi(v)\odot w$. In fact this is an element of
the corresponding tensor product.

Using Proposition 1 and its corollaries, it is easy to check
that the following proposition holds.

\begin{prop} The composition (product) of a 
\b-nuclear and a 
\b-linear mapping or of a \b-linear and a \b-nuclear mapping 
is a \b-nuclear operator.
\end{prop}

{\bf Proof.} It is obvious from our definitions that the composition
of a one-dimensional and a \b-linear mapping or of a 
\b-linear and a one-dimensional mapping 
is a one-dimensional mapping. So it is sufficient to
note that the decomposition with a \b-linear operator transforms each sum
of operators to a sum of operators. 
$\;\square$

\begin{trm} Suppose that $W$ is a \b-complete semimodule over a
\b-complete semiring $K$ and $V\subset K(X)$ is a functional \b-semimodule.
If every \b-linear functional on $V$ is integral, then any \b-linear
operator $A:V\to W$ has an integral representation if and only if
it is \b-nuclear.
\end{trm}

{\bf Proof .}
Note that all functionals of the form $\delta _x$ (for $x\in X$) are
\b-linear because $V$ is a functional \b-semimodule, see
Proposition 2 above. Each \b-linear functional is integral in
our case, so every \b-nuclear operator is a sum of a collection of
one-dimensional operators of the form
$k_w(x)\delta _x\odot w$, where $k_w\in K(X)$, $w\in W$. Therefore,
each \b-nuclear operator is of the form $\OSup _{x\in X, w\in W}k_w(x)
\delta _x\odot w$, where $k_w\in K(X)$, i.e., this operator has an
integral representation with the kernel $k(x)=\OSup_{w}k_w(x)\odot w$. 

On the other hand, each integral operator with a kernel $k:X\to W$ 
can be presented in the form $\OSup_{x\in X}A_x$, where the
one-dimensional operator $A_x$ is defined by the formula
 $A_x=\delta _x\odot k(x)$, so it is \b-nuclear.$\;\square$

\subsection{The b-approximation property and b-nuclear\\
 semimodules}
We shall say that a \b-complete semimodule $V$ has the {\it 
\b-approximation property} if the identity operator $id:V\to V$ is
\b-nuclear (for a treatment of the approximation property for
locally convex spaces in the traditional functional analysis see \cite{21, 22}).

Let $V$ be an arbitrary \b-complete semimodule over a \b-complete
idempotent semiring $K$. We call this semimodule a {\it \b-nuclear
semimodule} if any \b-linear mapping of $V$ to an arbitrary
\b-complete semimodule $W$ over $K$ is a \b-nuclear operator. Recall
that, in the traditional functional analysis, a locally convex space
is nuclear if and only if all continuous linear mappings of this
space to any Banach space are nuclear operators, see
\cite{21, 22}.

Using Propositions 1 and 5, Corollary 1, and the fact that
every mapping will not be changed after (left or right)
multiplication by the identity operator, from our basic
definitions it is easy to deduce the following proposition.

\begin{prop} Let $V$ be an arbitrary \b-complete semimodule over
a \b-complete semiring $K$. The following statements are equivalent:
\begin{itemize}
\item[{\rm (1)}] the semimodule $V$ has the \b-approximation property;
\item[{\rm (2)}] each b-linear mapping from $V$ to an arbitrary \b-complete
semimodule $W$ over $K$ is \b-nuclear;
\item[{\rm (3)}] each \b-linear mapping from an arbitrary \b-complete
semimodule $W$ over $K$ to the semimodule $V$ is \b-nuclear.
\end{itemize}
\end{prop}

{\bf Proof.} 
If the identity operator $id:V\to V$ is  \b-nuclear and $f:V\to W$ 
is an arbitrary \b-linear mapping, then from the equality
$f=f\circ id$ and Proposition 5 we deduce that the operator
$f$ is \b-nuclear. On the other hand, if each \b-linear mapping
from $V$ to an arbitrary \b-complete semimodule $W$ over $K$
is \b-nuclear, then this is true for the identity mapping  
$id:V\to V$. Therefore, statements (1) and (2) are equivalent.
The equivalence of statements (1) and (3) can be proved similarly.$\;\square$

\begin{cor} An arbitrary \b-complete semimodule over
a \b-complete semiring $K$ is \b-nuclear if and only if this
semimodule has the \b-approximation property.
\end{cor}

Recall that, in the traditional functional analysis, any
nuclear space has the approximation property but the converse
statement is not true.

Some concrete examples of \b-nuclear spaces and semimodules are
described in Examples 1, 2, 4, and 5 (see above in the end of
Section 1). Important \b-nuclear spaces and semimodules are 
described below in Section 3 (e.g., Lipschitz spaces and 
semi-Lipschitz semimodules over commutative semirings). It is
easy to show that the idempotent spaces $USC(X)$ and 
$Conc(\Bbb R)$  (see Examples 3 and 6) are not \b-nuclear
(however, for these spaces the kernel theorem is true). The reason 
is that these spaces are not functional \b-spaces and the
corresponding $\delta$-functionals are not \b-linear
(and even linear).

\subsection{Kernel theorems for functional b-semimodules}
Let $V\subset K(X)$ be a \b-complete functional semimodule
over a \b-complete semiring $K$. Recall that for $V$ the
{\it kernel theorem} holds if each \b-linear mapping of
this semimodule to an arbitrary \b-complete semimodule over
$K$ has an integral representation.

\begin{trm} Suppose that a  \b-complete semiring $K$ and a nonempty
set $X$ are given. The kernel theorem holds for any functional \b-semimodule
$V\subset K(X)$ if and only if each \b-linear functional on
$V$ is integral and the semimodule $V$ is \b-nuclear, 
i.e., it has the \b- approximation property.
\end{trm}

{\bf Proof.} The theorem follows from Theorem 2 and Proposition 6.
$\square$

\begin{cor} If for a functional \b-semimodule the kernel theorem
holds, then this semimodule is \b-nuclear.
\end{cor}

Note that the possibility to get an integral representation of
a functional means that it is possible to decompose it into a sum
of functionals of the form $\delta _x$.
The following example demonstrates that it is not always
possible to have an integral representation of a \b-linear
functional; moreover, this depends on embeddings of the 
semimodule to $K(X)$. On the other hand, the \b-approximation
property (\b-nuclearity) is invariant with respect to
isomorphism of semimodules.\\

{\bf Example 7.} Suppose that $K=\Rmax$, $X=\Bbb R$ and the function
$f\in K(X)$ is defined by the formula $f(x)=-x$, ?ਠ$x\in \Bbb R$. 
Denote by $V$ the subsemimodule of $K(X)$ consisting of all functions
of the form $a\odot f\oplus b$ for $a,b\in K$. It is easy to see
that $V$ is a \b-subsemimodule of $K(X)$ and the mapping
$(a,b)\mapsto a\odot f\oplus b$ is an isomorphism of the \b-space
$\Rmax\times \Rmax$ onto $V$. So the mapping $\phi$ transforming
$a\odot f\oplus b$ to $b$ is a \b-linear functional. Let us show
that $\phi$ has no integral representations. 

Indeed, let $k:X\to K$ be an integral kernel of the functional $\phi$. 
For an arbitrary $v=a\odot f\oplus\1\in V$ and $x\in X$ we have 
$0=\1=\phi (v)=\OSup_{y\in X} v(y)\odot k(y)\> v(x)\odot k(x)
= \max(a-x,0)+k(x)$, 
i.e., $0\>\max(a-x,0)+k(x)$, so  $k(x)\< -\max(a-x,0)$. The number 
$a$ can be chosen arbitrarily great, so $k(x)=\0$ for all $x\in X$, 
i.e. $k=\0$. But this is impossible because a nonzero functional
is not able to have the zero integral kernel. Therefore, the
functional $\phi$ has no 
integral kernels, as was to be proved.

{\bf Remark 3.} The semimodule presented in Example 7 is naturally
isomorphic to the semimodule $K(\{x,y\})$ of all functions defined on
the two-point set $\{x, y\}$ and the isomorphism $K(\{x,y\})\to V$ 
is defined by the formula $g\mapsto g(x)\odot f + g(y)$, where
$g\in K(\{x, y\})$.
So the functional $\phi$ described in Example 7 coincides with
$\delta_y$ and it is integral (with the kernel $\delta_y$) in
$K(\{x, y\})$. So we see that the property to be integral is not
an ``intrinsic'' property of a functional but depends on its 
imbedding to a semimodule of functions.\\

\noindent{\bf Theorem 3a} {\it Suppose that a \b-complete semiring $K$ 
and a nonempty set $X$ are given. The kernel theorem holds for a 
functional \b-semimodule
$V\subset K(X)$ if and only if the identity operator 
$id:V\to V$ is integral.}\\

{\bf Proof.} It follows from the obvious fact that the composition
(product) of any integral operator with each \b-linear
operator is an integral operator. Indeed, suppose that $A$ is a 
\b-linear operator transforming $V$ to a semimodule $W$ and 
$k:X\to V$ is an integral kernel of the operator $id$. Then 
$f=\OSup_{x\in X}f(x)\odot k(x)$ for each $f\in V$, so  
$(Af)(x)= \OSup_{x\in X}f(x)\odot A(k(x))$. Thus the mapping
$x\mapsto A(k(x))$ is a kernel of $A$. The converse statement
is trivial: the identity operator is integral if the kernel
theorem holds.$\;\square$

\section{A description of functional b-semimodules 
for which the kernel theorem holds}

Suppose that $X$ is a nonempty set and $K$ is an idempotent
semiring. We shall say that a function $d$ defined on
$X\times X$ and taking its values in $K$ is a 
{\it semimetric} on $X$ with values in $K$ if
\begin{equation}
d(x,y)=\OSup_{z\in X}d(x,z)\odot d(z,y), 
\end{equation}
where $x,y\in X$. We shall say that this semimetric is  
{\it symmetric} if $d(x,y) = d(y,x)$ for all $x, y \in X$.
If  $d(x,x)=\1$ for all $x\in X$, then the semimetric is
{\it reflexive}; in this case the condition (4) 
is equivalent to the
triangle inequality: $d(x,y)\> d(x,z)\odot d(z,y)$, where $x,y,z\in X$. 

{\bf Example 8.} Let $r=r(x,y)$ be a metric on $X$. Then the function 
$d(x,y)=-r(x,y)$ is a reflexive 
symmetric semimetric on $X$ with values in $\Rmax$.\\

Let us present an example of a nonsymmetric and nonreflexive
semimetric on the set of real numbers.\\

{\bf Example 9.} Suppose $X=\Bbb R$, $K=\Rmax$. Set $d(x,y)=\1$ if $y < x$ and
$d(x,y)=\0$ if $x\leq y$. Then $d(x,y)$ is a nonsymmetric semimetric and 
$d(x,x)=\0$ for all $x\in X$.\\

Now we introduce our basic definitions for this section. For
a semimetric $d$ on $X$ with values in a \b-complete semiring
$K$, we define \b-{\it closed functional semimodules} $\; \Lip(X,d)$ 
{\it and} $\; \lip(X,d)$ 
{\it on} $X$ by means of the formulas
$$
\begin{array}{rcl}
\Lip(X,d) & = &\{f\in K(X)|f(x)\> f(y)\odot d(y,x)\; \mbox{ for all } 
x, y\in X \},\\
\lip(X,d) & = &\{f\in K(X)|f(x)=\OSup_{y\in X}f(y) \odot d(y,x)\; 
\mbox{ for all } x\in X\}.\\ 
\end{array}
$$

The functional semimodules of the form $\lip(X,d)$ will be called
{\it semi-Lipschitz semimodules}.

It follows from the definition that $\lip(X,d)$ is a \b-closed
subsemimodule of $\Lip(X,d)$. It is easy to see that in the case
of a reflexive semimetric $\lip(X,d)$ and $\Lip(X,d)$ coincide. But
in general this is not true. In particular, in the situation  of Example
9 the space $\Lip(X,d)$ consists of all nonincreasing functions on
$\Rmax$ whereas $\lip(X,d)$ is the space of all lower semicontinuous
functions\footnote{Lower semicontinuous functions taking their
values in $\Rmax$ are defined in the same way as in Example 3 (see
Section 1.2 above), but $X_s$ is defined as the set
$\{x\in X\vert f(x)\leq s\}$; see also \cite{8}.}
belonging to $\Lip(X,d)$.\\

{\bf Example 10.} Let $X$ be a nonempty metric space with a
fixed metric $r$. Denote by $\lip(X)$ the set of all functions
defined on $X$, taking their values in $\Rmax$, and
satisfying the following {\it Lipschitz condition}:
$$
\mid f(x)\odot (f(y))^{-1}\mid \; = \;
\mid f(x) - f(y)\mid \; \leq \; r(x, y),
$$
where $x$, $y$ are arbitrary elements of $X$.  The set
$\lip(X)$ consists of continuous real-valued functions (but not all
of them!) and 
(by definition) the function which is equal to
$-\infty = \0$ at each point $x\in X$. It is easy to check that
$\lip(X) = \lip(X,d)$, where $d(x, y) = -r(x, y)$ (see Example 8
above); so $\lip(X)$ has a structure of an idempotent space over the
semifield $\Rmax$. Spaces of the form $\lip(X)$ are said to be
 {\it Lipschitz spaces} (see also \cite{8}, example 2.9.12). These
spaces are \b-nuclear and for each Lipschitz space the kernel
theorem holds (see Theorem 4 below).\\

For $x\in X$ denote by $d_x$ the function on $X$ defined by the
formula $d_x:y\mapsto d(x,y)$. By virtue of the equality (4)
all functions of the form $d_x$ are elements of the space $\lip(X,d)$. 
Denote by $\lip_0(X,d)$ the {\it subsemimodule} of
$\lip(X,d)$ {\it generated by these functions} (i.e. consisting
of all their finite linear combinations).

\begin{prop}
Suppose that a semimetric $d$ takes its values in a 
\b-complete semiring $K$. Each \b-subsemimodule $V$ of $\; \lip(X,d)$
such that $V\supset\lip_0(X,d)$ is a lower ideal in $\; \lip(X,d)$ in the
sense of the lattice theory {\rm \cite{33}}. This means that if a function $f$
is an element of $V$, then $V$ contains all the elements of
$\; \lip(X,d)$ majorized \footnote{Recall that an element $g$
is majorized by $f$ if $g\< f$.} by $f$.  In particular, $V$ is
closed in $\; \lip(X,d)$ 
(but not necessarily in $K(X)$)  under the operation of
taking infima (the greatest lower bounds) of arbitrary nonempty
subsets.
\end{prop}                                              

{\bf Proof.} Suppose $f\in V$, $g\in \lip(X,d)$ and $g\< f$. Then the
set of functions $S = \{g(x)\odot d_x|x\in X\}\subset \lip_0(X,d)\subset V$ 
is bounded in $V$ because all elements of this set are majorized
by $f$. Using our definition of $\lip(X,d)$,
for each $g\in \lip(X,d)$ we have $g(x) =  
\OSup_{y\in X} g(y)\odot d_y(x)$,
i.e. $g = \OSup_{y\in X} g(y)\odot d_y$ and $g = \OSup S$ is an
element of $V$.
$\;\square$

\begin{trm} Suppose that a \b-complete semiring
$K$ and a non\-e\-mpty set $X$ are given and $V$ is a nondegenerate
functional \b-subsemimodule of $K(X)$. For $V$ the kernel theorem
holds if and only if there exists a semimetric $d$ on $X$ such that
$V$ is a \b-subsemimodule of $\; \lip(X,d)$ and $V\supset \lip_0(X,d)$.
In particular, if $V$ is a \b-closed subsemimodule of $K(X)$,
then $V = \lip(X,d)$.
\end{trm}

{\bf Proof.} Suppose that the kernel theorem holds for $V$ and $d$
is a kernel of the identity operator $id: V\to V$. Then we can use
$d$ as a semimetric on $X$. Indeed, for each $x\in X$, the function
$d_x: y\mapsto d(x,y)$ is an element of $V$ by the kernel definition
and $d_x=id(d_x)=\OSup_{z\in X}d_x(z)\odot d_z$. This can be rewritten 
in the form 
$d(x, y) 
=\OSup_{z\in X}d(x,z)\odot d(z,y)$ for each $y\in X$
because $V$ is a functional \b-semimodule. So $d$ is a semimetric on $X$
and  $\lip_0(X,d)\subset V$. Then for each $f\in V$ we have 
$f=id(f)=\OSup_{z\in X}f(z)\odot d_z$; so 
$f(y)=\OSup_{z\in X}f(z)\odot d(z,y)$ for each $y\in X$. 
Thus $V\subset \lip(X,d)$.

On the other hand if $V$ is a \b-subsemimodule of $\lip(X,d)$ and $V$
contains $\lip_0(X,d)$, then each function of the form $d_x$ belongs
to $\lip_0(X,d)\subset V$. So, by our definition of $\lip(X,d)$, the 
corresponding semimetric $d$ is
a kernel of the identity mapping $id: V\to V$ and $d\in kern_V(X)$. $\;\square$

\begin{cor} Lipschitz spaces are nuclear. Nondegenerate semi-Lipschitz
semimodules over \b-complete semirings are nuclear.
\end{cor}

Note that Corollary 5 is also true for degenerate admissible
semi-Lipschitz semimodules over commutative semirings
(see Remark 5 below). \\

{\bf Remark 4.} If (under the conditions of Theorem 4) the \b-semimodule
 $V$ is admissible and it is a $\wedge$-semimodule, then
the semimetric $d$ constructed for the proof of Theorem 4 is
reflexive. In particular, if the \b-semimodule $V$ is nondegenerate,
then the maximal kernel of the identity operator (this integral
kernel is constructed for the proof of Theorem 4) is a reflexive semimetric.\\

{\bf Remark 5.} Suppose that a functional \b-semimodule $V\subset K(X)$ 
is degenerate but admissible (for example, $V$ is admissible
automatically if $K$ is a semifield). Then an analog of the
Theorem 4 is true for the semimodule consisting of all
restrictions from $V$ to the set 
$X_V=\{ x\in X|(\exists f\in V){ : }f(x)  = \1\}$, see Proposition 3.

\section{Integral representations of operators in\\
 abstract idempotent semimodules}

In this section we examine the following problem: when
a \b-complete idempotent semimodule $V$ over a \b-complete
semiring is isomorphic to a functional \b-semimodule
$W$ such that the kernel theorem holds for $W$.

Suppose that $V$ is a \b-complete idempotent
semimodule over a \b-complete semiring $K$ and $\phi$ is a 
\b-linear functional defined on $V$. We call this
functional a $\delta$-{\it functional} if there exists an 
element $v\in V$ such that 
$$
\phi(w)\odot v\< w
$$ 
for each element $w\in V$. It is easy to see that every
functional of the form $\delta_x$ is a $\delta$-functional
in this sense (but the converse is not true in general).

Denote by $\Delta(V)$ the set of all $\delta$-functionals on 
$V$. Denote by $i_\Delta$ the natural mapping
$V\to K(\Delta(V))$defined by the formula
$$
(i_\Delta (v) )(\phi)=\phi(v) 
$$
for all $\phi \in \Delta(V)$.
We shall call an element $v\in V$ {\it pointlike} if there
exists a \b-linear functional $\phi$ such that 
$\phi(w)\odot v\< w$ for all $w\in V$. The set of all
pointlike elements of $V$ will be denoted by $P(V)$. Recall
that by $\phi\odot v$ we denote the one-dimensional operator
$w\mapsto \phi(w)\odot v$. 

The following statement is an
obvious consequence of our definitions (including the
definition of the standard order) and idempotency
of our addition. \\

{\bf Remark 6.} If a one-dimensional operator $\phi\odot v$ 
appears in a decomposition of the identity operator on $V$
into a sum of one-dimensional operators, then $\phi\in\Delta(V)$ 
and $v\in P(V)$.\\

Denote by $id$ and $Id$ the identity operators on $V$ and 
$i_\Delta (V)$, respectively.

\begin{prop} {\rm 1)} If the operator $id$ is \b-nuclear, then
$i_\Delta$ is an embedding and the operator $Id$ is integral.

{\rm 2)} If the operator  $i_\Delta$ is an embedding and the operator $Id$ 
is integral, then the operator $id$ is \b-nuclear.
\end{prop}

{\bf Proof.} Statement 2) is obvious, so it is sufficient to prove
statement 1). Then (by the condition of this statement) the
operator $id$ is a sum of one-dimensional operators, i.e. 
$id = \OSup_{\phi\in \Delta}\phi\odot w_\phi$ for some collection 
of elements $w_\phi \in V$. In this case we have 
$Id = \OSup_{\phi\in \Delta}\delta_\phi\odot i_\Delta (w_\phi)$, i.e. 
the operator $Id$ is integral, as was to be proved. $\square$

The following result is a direct consequence of Proposition
8, Theorem 3a, and Proposition 6.

\begin{trm} A \b-complete idempotent semimodule $V$ over
a \b-complete idempotent semiring $K$ is isomorphic to a
functional \b-semimodule for which the kernel theorem holds
if and only if the identity mapping on $V$ is a
\b-nuclear operator, i.e. $V$ is a \b-nuclear semimodule.
\end{trm}

The following proposition shows that, in a certain sense, the
imbedding $i_\Delta$ is a universal representation of a 
\b-nuclear semimodule in the form of a functional
\b-semimodule for which the kernel theorem holds. 

\begin{prop} Let $K$ be a \b-complete idempotent semiring, $X$ 
a nonempty set, and $V\subset K(X)$ a functional \b-semimodule
on $X$ for which the kernel theorem holds. Then there exists
a natural mapping $i:X\to \Delta(V)$ such that the corresponding
mapping $i_*: K(\Delta(V))\to K(X)$ is an isomorphism of $i_\Delta(V)$
onto $V$.
\end{prop}

{\bf Proof.} The imbedding $i:X\to \Delta(V)$ is defined by the
formula $i: x\mapsto \delta_x$, $x\in X$. This definition is correct by
virtue of Proposition 1. It is easy to check that the mapping 
$i_*: K(\Delta(V))\to K(X)$ is an isomorphism of $i_\Delta(V)$ onto $V$. 
$\;\square$

{\bf Remark 7.} For the sake of simplicity, we treat 
only \b-complete semimodules in this paper. Using the procedure
of bounded completion and operating in the spirit of the paper \cite{8},
it is possible to extend a considerable part of the
definitions and results to the case of incomplete semimodules
over incomplete semirings.

\medskip

The authors are sincerely grateful to V.N. Kolokoltsov,
V.P. Maslov, A.N. Sobolevski{{\u\i}}, and A.M. Vershik for valuable
suggestions and support.

\end{document}